\documentclass{amsart}

\usepackage{amsmath,amssymb,amsfonts}

\usepackage{mathtools}
\usepackage[shortlabels]{enumitem}
\usepackage{xcolor}
\usepackage{mathrsfs}
\usepackage{nicefrac}
\usepackage{Orcidlinknew}

\usepackage{natbib}
 \bibpunct[, ]{(}{)}{,}{a}{}{,}%
 \def\BIBand{and}%

\usepackage{rotating}
\usepackage{fancyvrb}

\definecolor{linkcolor}{RGB}{119,30,17}

\usepackage{hyperref}

\definecolor{pnas}{RGB}{185,69,35}
\definecolor{pnasblue}{RGB}{12,100,180}
\definecolor{el}{RGB}{147,72,29}
\definecolor{elo}{RGB}{191,90,35}

\hypersetup{
    colorlinks=true,
    linkcolor=el,
    filecolor=el,      
    urlcolor=el,
    citecolor=el,
    }

 \DeclareMathSymbol{\Nset}{\mathbin}{AMSb}{"4E}
 \DeclareMathSymbol{\Zset}{\mathbin}{AMSb}{"5A}
 \DeclareMathSymbol{\Rset}{\mathbin}{AMSb}{"52}
 \DeclareMathSymbol{\Qset}{\mathbin}{AMSb}{"51}
  \DeclareMathSymbol{\Fset}{\mathbin}{AMSb}{"46}
 \DeclareMathSymbol{\Cset}{\mathbin}{AMSb}{"43}
 \DeclareMathSymbol{\Kset}{\mathbin}{AMSb}{"4B}
 
 \DeclareMathSymbol{\Sset}{\mathbin}{AMSb}{"53}

\newcommand\ind[1]{\ensuremath{\mathbf{1}_{#1}}}
\newcommand\dm{\ensuremath{\rho}}
\newcommand{\term}[1]{\textbf{#1}}

\definecolor{defcolor}{RGB}{14,45,97}
\definecolor{thmcolor}{RGB}{119,30,17}
\newcommand{\dflist}{\protect\color{defcolor}} 
\newcommand{\thmlist}{\protect\color{thmcolor}} 

\newcommand\xqed[1]{%
  \leavevmode\unskip\penalty9999 \hbox{}\nobreak\hfill
  \quad\hbox{{#1}}}
\newcommand\qeddef{\xqed{${\scriptstyle\spadesuit}$}}
\newcommand\qedthm{\xqed{${\scriptstyle\blacksquare}$}}
\newcommand\qedexa{\xqed{$\scriptstyle\blacklozenge$}}

\theoremstyle{definition}
\newtheorem{theorem}{Theorem}

\newtheorem{example}[theorem]{Example}

\newtheorem{definition}{Definition}

\newcommand\iless{\blacktriangleleft}
\newcommand\ccon[2][\ensuremath{>}]{\textcolor[RGB]{14,45,97}{\sc \selectfont \ensuremath{#1}C{#2}}}

   \newcommand{\dom}[2][\big]{\ensuremath{\mathbin{\mathrm{dom}{\mathopen{#1(\,}\mathop{#2}\mathclose{\,#1)}}}}}
    \newcommand{\img}[2][\big]{\ensuremath{\mathbin{\mathrm{im}{\mathopen{#1(\,}\mathop{#2}\mathclose{\,#1)}}}}}

 \DeclareMathOperator{\sgn}{sgn}  

   \newcommand{\opposite}[2][]{\ensuremath{\mathbin{\mathrm{op}{\mathopen{#1(\,}\mathop{#2}\mathclose{\,#1)}}}}}
   \newcommand{\ap}[2][]{\ensuremath{\mathbin{\mathrm{ap}{\mathopen{#1(\,}\mathop{#2}\mathclose{\,#1)}}}}}

\setlist[itemize]{leftmargin=*} 
\newlist{myitemize}{itemize}{3}
\setlist[myitemize,1]{label=\textbullet,leftmargin=*}
\setlist[myitemize,2]{label=$\rightarrow$,leftmargin=4em}
\setlist[myitemize,3]{label=$\diamond$}

\makeatletter
\newcommand{\vast}{\bBigg@{3}}
\newcommand{\Vast}{\bBigg@{4}}
\newcommand{\VVast}{\bBigg@{5}}
\newcommand{\vastl}{\mathopen\vast}

\newcommand{\vastr}{\mathclose\vast}

\makeatother

\begin{document}



\title[Formal Power Series Representations in Probability]{Formal Power Series Representations in Probability and Expected Utility Theory}

%
\author{\orcidlinki{Arthur Paul Pedersen}{0000-0002-2164-6404}$^\dagger$ \and \orcidlinki{Samuel Allen Alexander}{0000-0002-7930-110X}$^\ddagger$}
\email[Arthur Paul Pedersen]{appedersen@cs.ccny.cuny.edu}

\email[Samuel Allen Alexander]{samuel.allen.alexander@gmail.com}




\noindent\thanks{$^\dagger$Department of Computer Science, The City University of New York.\\
\hspace*{12pt}$^\ddagger$Quantitative Analytics Unit, U.S. Securities \& Exchange Commission}






\begin{abstract}We advance a general theory of coherent preference that surrenders restrictions embodied in orthodox doctrine. This theory enjoys the property that any preference system admits extension to a complete system of preferences, provided it satisfies a certain coherence requirement analogous to the one de Finetti advanced for his foundations of probability.  Unlike de Finetti's theory, the one we set forth requires neither transitivity nor Archimedeanness nor boundedness nor continuity of preference. This theory also enjoys the property that any complete preference system meeting the standard of coherence can be represented by utility in an ordered field extension of the reals. Representability by utility is a corollary of this paper's central result, which at once extends Hölder's Theorem and strengthens Hahn's Embedding Theorem.
\end{abstract}

\keywords{comparative probability; qualitative probability; decision theory; utility theory; Bruno de Finetti; coherence; non-Archimedean systems; formal power series;  incomplete preference; indeterminacy; expected utility; H\"{o}lder's Theorem; Hahn's Embedding Theorem; representation theorem}


\maketitle

\section{Introduction}\label{sec1}

This paper shows that any complete system of preferences 
admits a utility representation in an ordered field
extension of the real numbers, subject to a certain standard of coherence analogous to one for gambling that \citet{deFinetti:1931a,deFinetti:1937,deFinetti:1974a,deFinetti:1974b} advanced for his theory of probability.   In so doing, this paper provides a full, self-contained proof of a principal result due to one of us but for space considerations could only be stated without proof in a prior publication \citep{Pedersen:2014}.  More generally,  this paper provides for a unified treatment of numerical probability and expected utility that relaxes standards imposed by the orthodox canon to which belong the staple theories of \citet{deFinetti:1974a,deFinetti:1974b},  \citet{vonNeumannMorgenstern:1947}, 
 \citet{Savage:1954}, and  \citet{AnscombeAumman:1963}, while making room for adherence to standards whose transgression is all but prescribed by the orthodox canon.\footnote{Among orthodox standards at issue are continuity requirements, Archimedean postulates, topological conditions, boundedness constraints, cardinality restrictions, and the like. Relaxing such requirements makes room for adherence to criteria requiring respect for forms of monotonicity or dominance preservation.  For more on the challenges the unified theory addresses, see exposition 
by \citet{PedersenArlo:2012,PedersenDissertation:2013,Pedersen:2014}. An abridged technical orientation to the result stated but unproven by \citet{Pedersen:2014} will be given  given in \S\ref{sec:background}. See footnote \ref{fn:specialissue}.}

 
 For the purpose of illustration by comparatively direct and elementary  means, this paper focuses on de Finetti's foundations of probability and expected utility, thereby sidestepping a detour through, for example, mixture spaces, convex sets, cones, affine spaces, and so on, in order to  formulate and apply this paper's result to extending the theories of \citet{deFinetti:1974a,deFinetti:1974b}, \citet{vonNeumannMorgenstern:1947}, \citet{AnscombeAumman:1963}, and
 \citet{Savage:1954}. Routine arguments show how to used to to in this connection, see, e.g.,  \citet{StoneM.H.1949Pftb}, \citet{Hausner:1954}, \citet{GudderStanleyP.1977CaM,GudderStanleyP.1978ECaM}, \citet{GudderS.1980GC}, \citet{MonginPhilippe2001Anom}).

The unified treatment of numerical probability and expected utility representations is a basic application of this paper's central result. In what follows, we show that any totally ordered linear space $\mathbb{V}$ is order isomorphic to a subspace of a lexicographically ordered field $\Rset\bigl(\bigl(\,\epsilon^{\Gamma}\,\bigr)\bigr)$  of formal power series $s=\sum_{\gamma}s_{\gamma}\epsilon^{\gamma}$ with coefficients $s_{\gamma}$ in $\Rset$ and exponents $\gamma$ in a totally ordered Abelian group $\Gamma$  for which $(\gamma:s_{\gamma}\neq 0)$ forms a well-ordered subset of $\Gamma$.  The operations of addition and multiplication agree with the familiar operations for power series subject to the usual requirement that  $\epsilon^{\alpha}\epsilon^{\beta}=\epsilon^{\alpha+\beta}$ for $\alpha,\beta\in \Gamma$, while the lexicographic order on $\Rset\bigl(\bigl(\,\epsilon^{\Gamma}\,\bigr)\bigr)$ is determined by requiring a power series $s$ to be positive when the first exponent $\gamma$ at which the real coefficient $s_{\gamma}$ does not vanish to be positive.  So equipped, the field $\Rset\bigl(\bigl(\,\epsilon^{\Gamma}\,\bigr)\bigr)$ is a totally ordered extension of the real number system. 

When $\Gamma$ is the trivial group (i.e., $\Gamma=\{0\}$), the field $\Rset\bigl(\bigl(\,\epsilon^{\Gamma}\,\bigr)\bigr)$
 is isomorphic to the field $\Rset$, in which case both $\mathbb{V}$ and $\Rset\bigl(\bigl(\,\epsilon^{\Gamma}\,\bigr)\bigr)$ are \emph{Archimedean} in the sense that  each pair $(a_{1},a_{2})$  of non-zero elements satisfies $n|a_{1}|>|a_{2}|$ for some $n\in\Nset$,
where $|a|\coloneqq \max\{a, -a\}$. Thus, a corollary of this paper's central result is H\"{o}lder's Theorem \citeyearpar{Holder:1901}, a cornerstone of modern measurement theory. When the totally ordered Abelian group $\Gamma$ is non-trivial (i.e., $\Gamma\neq \{0\}$), the field $\Rset\bigl(\bigl(\,\epsilon^{\Gamma}\,\bigr)\bigr)$ is a proper extension of the field of real numbers, in which case both $\mathbb{V}$ and $\Rset\bigl(\bigl(\,\epsilon^{\Gamma}\,\bigr)\bigr)$ are \emph{non-Archimedean}.  This paper's central result therefore strengthens Hahn’s Embedding Theorem \citeyearpar{Hahn:1927}. In fact, the field $\Rset\bigl(\bigl(\,\epsilon^{\Gamma}\,\bigr)\bigr)$ is, up to order-isomorphism, the smallest totally ordered field extension $\Fset$ of $\Rset$ to include (an order-isomorphic copy of) the linear space $\mathbb{V}$ that is \emph{Archimedean complete} in the sense that it admits no  proper extension $\Fset^{\prime}$ for which each element $a'\in \Fset^{\prime}\setminus\Fset$ has the same order of magnitude as some  $a \in \Fset$ --- that is, such that $n|a|>|a'|$ and $n|a'|>|a|$ for some $n\in\Nset$. 

The rest of this paper is organized as follows. After briefly laying out  preliminary terminology and notation in \S\ref{sec:prelims}, attention is turned in \S\ref{sec:background} to this paper's fundamental results on coherence and discusses their connection to pioneering work of de Finetti. 
Examples given in \S\ref{sec:examples} shed light on the scope and power of our theory of coherence, whereupon  \S\ref{sec:Representability} develops critical ideas and conventions for this paper's central result.  Concluding the paper in \S\ref{sec:remarks} is a brief discussion of related literature. 

\section{Preliminaries}
\label{sec:prelims}
Throughout we will be concerned with real-valued maps on a common nonempty set $\Omega$ of \term{states} --- called a \term{state space} --- corresponding to a collection of mutually exclusive and collectively exhaustive hypotheses. A function $g:\Omega\to\Rset$, called a \term{gamble}, specifies the numerical outcome to be obtained from $g$ in any given state.  A gamble only assuming the values $0$ or $1$ is called an \term{indicator function}. Each subset $E$ of $\Omega$, called an \term{event}, corresponds to a unique indicator function, denoted by $\ind{E}$, such that for every
$\omega\in\Omega$:
        \begin{align*}
\ind{E}(\omega)\quad&= \quad \begin{cases}
1 & \mbox{if }\omega\in E;\\
0 &\mbox{otherwise.}
\end{cases}
\end{align*}
Each indicator function likewise corresponds to a unique event from $\Omega$. Following ordinary convention,  an event $E$ will be identified with its indicator function $\ind{E}$ when the context is clear.  

A gamble $g$ is said to be \term{bounded} if its image is a subset of $\bigl[\,-n,n\,\bigr]$ for some nonnegative integer $n$; a gamble that fails to be bounded is said to be \term{unbounded}.  Observe that any gamble assuming only finitely many values is bounded; thus any linear combination of indicators functions is bounded.  Observe in addition that a gamble is unbounded only if the underlying state space is infinite.

In this paper, we presume gamble outcomes are expressed in units of a  linear utility scale that has been determined in advance. We also presume individual gambles can be combined and rescaled in accordance with pointwise arithmetic operations.  Thus,  given gambles $f$ and $g$ on state space $\Omega$ and $c\in \Rset$:
    \begin{enumerate}[\dflist (i)]
   \item
        The \term{sum} $f+g$, \term{difference} $f-g$, and \term{product} $fg$ of $f$ and $g$ on $\Omega$ are gambles on $\Omega$ such that for every $\omega\in\Omega$:
        \begin{align*}
\Bigl(\,f+g\,\Bigr)(\omega)\quad&\coloneqq \quad f(\omega)\,+\,g(\omega);\\
\Bigl(\,f-g\,\Bigr)(\omega)\quad&\coloneqq \quad f(\omega)\,-\,g(\omega);\\
\Bigl(\,fg\,\Bigr)(\omega)\quad&\coloneqq \quad f(\omega)g(\omega).
        \end{align*}
           \item
        The \term{scalar product} of $c$ and $f$ on $\Omega$ is the gamble $cf$ on $\Omega$ such that for every $\omega\in\Omega$:
        \begin{align*}
\Bigl(\,cf\,\Bigr)(\omega)\quad&\coloneqq \quad cf(\omega).
        \end{align*}
          \item The \term{constant} $c$ \term{gamble} on $\Omega$ is the gamble $\mathbf{c}$ on $\Omega$ such that $\mathbf{c}(\omega)\coloneqq c$ for all $\omega\in\Omega$.
    \end{enumerate} 
Observe that the constant $1$ gamble $\mathbf{1}$ is the indicator function $\ind{\Omega}$ for $\Omega$, the \term{sure event}, while the constant $0$ gamble $\mathbf{0}$ is the indicator  function $\ind{\varnothing}$ for $\varnothing$, the \term{impossible event}.
    
To fix ideas, consider state space $\Omega=\Bigl\{\,\mathtt{Heads},\mathtt{Tails}\,\Bigr\}$
and gamble $h$ with $h\bigl(\mathtt{Heads}\bigr)\,=\, -h\bigl(\mathtt{Tails}\bigr)\,=\,1$, representing the uncertain reward contingent on the result of a coin toss:  If the coin lands heads, then the gamble brings in
 $\$1$; but if it lands tails, then the gamble debits $\$1$. Observe that gamble $h$ is statewise better than the constant gamble $-\mathbf{c}$  for any real number $c>1$ that gamble $h$ results in a higher reward than the constant gamble $-\mathbf{1}$ if the coin lands heads but results in the same reward as the constant gamble $-\mathbf{1}$ if the coin lands tails. We introduce terminology for such statewise relationships among gambles in the following definition.
 
\begin{definition}[Dominance]
    \label{gamblebasicsrelsdefn}
    Let $f$ and $g$ be gambles on state space $\Omega$.
    \begin{enumerate}[\dflist (i), leftmargin=0.15in,labelsep=0.25in,itemindent=0.5in]
         
        \item
         Gamble $f$ is said to \term{dominate} gamble  $g$, abbreviated $g\leq f$, if for all $\omega\in\Omega$, $g(\omega)\leq f(\omega)$.
         \item    Gamble $f$ is said to \term{uniformly dominate} gamble  $g$, abbreviated $g\ll f$, if $g\,+\mathbf{c}\leq f$ for some positive real number $c$.
\item         Gamble $f$ is said to \term{weakly dominate} gamble  $g$, abbreviated $g< f$, if $g(\omega)\leq f(\omega)$ for all $\omega\in\Omega$ and $g(\omega)< f(\omega)$ for some $\omega\in\Omega$.
         \item   Gamble $f$ is said to \term{simply dominate} gamble  $g$, abbreviated $g \lessdot  f$, if for all $\omega\in\Omega$, $g(\omega)< f(\omega)$.

          \qeddef{}
    \end{enumerate}
\end{definition}
Weak dominance of $f$ over $g$ requires $g\leq f$ and $g\neq f$. In the preceding example, gamble $h$ dominates, uniformly dominates, weakly dominates, and simply dominates  the constant gamble $-\mathbf{c}$  for any real number $c>1$; however, while gamble $h$ both dominates and weakly dominates  constant gamble $-\mathbf{1}$ , it neither uniformly dominates nor simply dominates constant gamble $-\mathbf{1}$.

Recall that a pair $(X,R)$ of sets for which $R\subseteq X\times X$ is called a \term{binary relation} on $X$.  A binary relation $(X,R)$ is said to be \term{complete} (or \term{total} or \term{linear}) if all elements of $X$ are $R$-comparable --- that is to say, if for all $a,b\in X$, either $(a,b)\in R$ or $(b,a)\in R$; called a (\term{strict preorder}) \term{preorder} (or \term{quasiorder}) if it is  (irreflexive and transitive) reflexive and transitive ;called a \term{weak order} if it is a complete preorder; called a \term{partial order} if it is a antisymmetric preorder; called a \term{total} (or \term{linear} or \term{complete}) \term{order}, or simply an \term{order}, if it is an total partial order. 

The \term{opposite} (or \term{dual} or \term{inverse} or \term{transpose}) of a binary relation $R$, denoted by $\opposite{R}$,  is the binary relation $\opposite{R}\coloneqq\bigl\{\,(a,b)\,:\,(b,a)\in R\,\bigr\}$, while the \term{asymmetric part} (or \term{asymmetric component}) of $R$, denoted by $\ap{R}$, is the binary relation $\ap{R}\coloneqq R\setminus \opposite[\Big]{R}$.

A (\term{totally ordered}, \term{partially ordered}, \term{totally preordered}, \term{strict preordered}) \term{preordered linear space} is a (total order, partial order, total preorder, strict preorder) preorder $(V,R)$ on a linear (or vector) space $V$ whose additive and scalar multiplicative operations are compatible with $R$ in the sense that for all $x,y,z\in V$ and positive scalars $\alpha\in\Rset$, $(x,y)\in R$ if and only if $\bigl(\,\alpha x+z,\alpha y+z\,\bigr)\in R$. 

A binary relation $(X_{2}, R_{2})$ is said to \term{extend} a binary relation $(X_{1}, R_{1})$ --- and the binary relation $(X_{2}, R_{2})$ is thereby said to be an \term{extension} of the binary relation $(X_{1}, R_{1})$ --- if $X_{1}\subseteq X_{2}$, $R_{1}\subseteq R_{2}$, and $\ap{R_{1}}\subseteq \ap{R_{2}}$. When $(X_{1}, R_{1})$ and $(X_{2}, R_{2})$ are irreflexive over a common domain $X=X_{1}=X_{2}$, binary relation $R_{2}$ is said to \term{respect}  binary relation $R_{1}$ if $(X_{2}, R_{2})$ extends $(X_{1}, R_{1})$.

We follow the usual convention of identifying a binary relation $(X,R)$ with its second coordinate $R$ when so abusing terminology presents no danger of confusion. Thus a subset $R$ of the Cartesian product is called transitive (reflexive, a weak order, etc.) if the binary relation $(X, R)$  is transitive (reflexive, a weak order, etc.).   We  also adopt the standard convention of using infix notation, writing $aRb$ to mean that $(a,b)\in R$.

 If $\varphi(x)$ expresses a property of elements $x$ from $X$, we also follow the usual convention of writing $\bigl(\,\varphi\,\bigr)$ as shorthand for $\{ x\in X: \varphi(x) \}$ when there is no danger of confusion.  Thus, for example, if $f$ and $g$ are gambles on state space $\Omega$ and $c\in \Rset$, then $\Bigl(\,f\neq g\,\Bigr) = \Bigl\{\,\omega\in\Omega:f(\omega)\neq g(\omega)\,\Bigr\}$ and $\Bigl(\,f\,<\,c\,\Bigr) = \Bigl\{\,\omega\in\Omega:f(\omega)\,< c\,\Bigr\}.$

\section{Background}
\label{sec:background}

We presume that an individual subscribes to a comparative preference ranking over gambles.

\begin{definition}[Comparative Expectations]
    Let $\mathscr X$ be a set of gambles on state space $\Omega$.
     A \term{comparative expectation system} on $\mathscr X$ is a pair of the form $\bigl(\mathscr{X},\succsim\bigr)$ such that $\succsim$ is a binary relation on $\mathscr X$.\qeddef{}
\end{definition}
The relation $f\succsim g$ registers the individual's \emph{weak preference} of $f$ over $g$. Operationally, the relation $f\succsim g$ is understood to express the individual's commitment to honoring a contract to exchange gamble $g$ in order to receive gamble $f$. The combination $f\succsim g$ and $g\succsim f$, abbreviated $f\sim g$, is understood to express the individual's \emph{indifference} between $f$ and $g$.

De Finetti \citeyearpar{deFinetti:1974a} requires each gamble $f$ from $\mathscr{X}$ to be associated with a unique real number $p$ such that $f\sim \boldsymbol{p}$. The value $p$ associated with $f$, denoted by $\mathbb{P}(f)$, is the \term{price} at which the individual is prepared to exchange gamble $f$. Thus, the value $\mathbb{P}(f)$ is the rate the individual judges exchanging any number of units of $f$ to be \emph{fair}.  De Finetti presumes that the individual's preference ranking over gambles accords with the ordering determined by their fair prices --- that is, $\mathbb{P}$ \emph{represents} $\succsim$ on $\mathscr{X}$ in the sense that for each $f,g\in\mathscr{X}$: $f\succsim g$ if and only $\mathbb{P}(f)\geq \mathbb{P}(g)$. Hence, according to de Finetti's treatment, all gambles from $\mathscr{X}$ are $\succsim$-comparable in the sense that either  $f\succsim g$ or  $g\succsim f$ for all $f,g\in\mathscr{X}$. When all gambles from $\mathscr{X}$  are $\succsim$-comparable, the comparative expectation system $\bigl(\,\mathscr{X},\succsim\,\bigr)$ is said to be \term{complete}.

\begin{definition} \label{ComparativePrevisionSystemdefn} A \term{comparative prevision system} is a comparative expectation system of the form $\bigl(\mathscr{X},\succsim\bigr)$ such that:  
\begin{enumerate}[\dflist \scshape dF1, leftmargin=0.2in,labelsep=0.25in,itemindent=0.5in,topsep=.08in,itemsep=0.1in]
\item Every gamble from $\mathscr{X}$ is bounded;
\item $\bigl(\mathscr{X},\succsim\bigr)$ is complete;
\item For each gamble $f\in\mathscr{X}$ there is $p\in \Rset$ with $\boldsymbol{p}\in\mathscr{X}$ such that $f\sim \boldsymbol{p}$.\qeddef{}
    \end{enumerate}
\end{definition}

De Finetti's criterion of coherence requires preferences to be secure against entering into a contractual arrangement resulting in a uniform net loss.  The following definition subsumes de Finetti's requirement as a special case.

\begin{definition}[$\dm$-Coherence] \label{coherentComparativeExpectationDefn} Let $\dm$ be a strict linear preorder respected by $>$ on the set of all gambles. A comparative expectation system $\bigl(\mathscr{X},\succsim\bigr)$ is said to be $\dm$-\term{coherent} if no positive integer $n$ exists for which
   there are gambles $f_1,\ldots,f_n,
    g_1,\ldots,g_n\in\mathscr X$ and positive scalars $c_1,\ldots,c_n\in\Rset$ such that:
    \begin{myitemize}[leftmargin=0.5in,labelsep=0.25in,itemindent=0.5in,itemsep=0.19in,topsep=.15in]
           \item[{\textcolor[RGB]{14,45,97}{\textsc{$\dm$C1}}}] \quad\quad $f_i\,\succsim\, g_i$\qquad\; for each $i=1,\ldots,n;$
       \item[{\textcolor[RGB]{14,45,97}{\textsc{$\dm$C2}}}] $\displaystyle\sum_{i=1}^n c_ig_i\, \geq\, \sum_{i=1}^n c_if_i$;
         \item[{\textcolor[RGB]{14,45,97}{\textsc{$\dm$C3$\mid$4}}}] At least one of the following two conditions obtain:
        \begin{myitemize}[leftmargin=0.5in,labelsep=0.25in,itemindent=0.5in,topsep=.15in,itemsep=0.19in]
              \item[{\textcolor[RGB]{14,45,97}{\textsc{$\dm$C3}}}] 
                \quad\quad $g_k\,\not\succsim\, f_k$ \qquad\; for some $k=1,\ldots,n;$
           \item[{\textcolor[RGB]{14,45,97}{\textsc{$\dm$C4}}}] $\displaystyle\sum_{i=1}^n c_ig_i \;\dm\; \sum_{i=1}^n c_i f_i$.
        \end{myitemize}
    \end{myitemize}
    The comparative expectation system $\bigl(\mathscr{X},\succsim\bigr)$ is said to be $\dm$-\term{incoherent}
    if it fails to be $\dm$-coherent.  We will omit the prefix $\dm$ when the intended relation is clearly understood from context. \qeddef{}
\end{definition}

Call a comparative prevision system \term{de Finetti coherent} if it is $\gg$-coherent.   We now show that the fundamental results due to de Finetti can be obtained when recast in the present framework:
\begin{enumerate}[\thmlist 1.,leftmargin=0.45in,labelsep=0.25in,itemindent=0.0in,itemsep=0.15in,topsep=.05in]
 \item If $f$ is a bounded gamble, then a de Finetti coherent comparative prevision system $\bigl(\mathscr{X},\succsim\bigr)$ can be extended to a de Finetti coherent comparative prevision system $\bigl(\,\mathscr{X}\cup\{f, \boldsymbol{p}_{f}\} ,\succsim^{\prime}\,\bigr)$ with $f\sim^{\prime} \boldsymbol{p}_{f}$, and the admissible values for ${p}_{f}$ can be characterized \cite[pp. 336-338]{deFinetti:1974b} (cf. \citet[pp. 105-107]{deFinetti:1949}, \citep[pp. 111-116]{deFinetti:1974a}).
 \item The price function $\mathbb{P}$ associated with a de Finetti coherent comparative prevision system  enjoys the properties of a \term{real-valued finitely additive expectation} \cite[pp. 69ff.]{deFinetti:1974a} (cf. \cite[pp. 105-106]{deFinetti:1949}, \cite[p. 335]{deFinetti:1974b}):    
\begin{figure}[h!]
\noindent
    \begin{tabular}{c c}
  \begin{minipage}{0.51\textwidth}\begin{enumerate}[\dflist \textsc{P}1,leftmargin=0.25in,labelsep=0.25in,itemindent=0.25in,itemsep=0.15in]
    \item $\mathbb{P}(\mathbf{1})\;=\;1$
    \item $\mathbb{P} \bigl(\,rf\,+\,sg\,\bigr)\;=\;r\mathbb P(f)\,+\,s\mathbb P(g)$
    \item If $f\gg g$, then
            $\mathbb{P}(f)\,>\,\mathbb{P}(g)$
  \end{enumerate}\end{minipage} &
  \begin{minipage}{0.5\textwidth}\begin{enumerate}[label={\strut},itemsep=0.15in, leftmargin=.2in]
    \item if $\mathbf{1}\in\mathscr{X}$;
    \item for all
            $r,s\in\Rset,f,g,rf+sg\in \mathscr{X}$;
    \item for all $f,g\in \mathscr{X}$.
  \end{enumerate}\end{minipage} 
\end{tabular}
\end{figure}
\end{enumerate}

The first result is a reformulation of de Finetti's \emph{Fundamental Theorem of Prevision}. Assuming without loss of generality that $\mathscr{X}$ includes all constant gambles, a corollary of this one-off extension result is that any de Finetti coherent comparative prevision system on $\mathscr{X}$ can be extended to a de Finetti coherent comparative prevision system on any  family of bounded gambles $\mathscr{Y}\supseteq \mathscr{X}$. 

The second result can be strengthened into a characterization of de Finetti coherence under an additional proviso: If $\mathscr{X}$ is a linear space containing all constant gambles, then a comparative prevision system $(\mathscr{X},\succsim)$ is de Finetti coherent if and only if the associated price function $\mathbb{P}$ is a real-valued finitely additive expectation that represents $\succsim$.

Since the indicator function for any event is a bounded gamble, it follows from the second result that the price function associated with a de Finetti coherent comparative prevision system satisfies the elementary properties of a real-valued finitely additive probability function when restricted to indicator functions.  

We relax the requirements of a real-valued finitely additive expectations to accommodate coherent preference.

\begin{definition}[Finitely Additive Expectation] \label{bbEdefn} Let $\Fset$ be a totally ordered field extension of  $\Rset$, and let $\mathscr X$ be a collection of gambles.
   An $\Fset$-\term{valued expectation} on $\mathscr{X}$ is a function $\mathbb{E}:\mathscr{X}\to\Fset$ such that:

\begin{enumerate}[\dflist \textsc{e}1,leftmargin=0.2in,labelsep=0.25in,itemindent=0.5in,topsep=.08in,itemsep=0.1in]
 \item $\mathbb{E}(\mathbf{1})\;=\;1$
    \item $\mathbb{E} \bigl(\,rf\,+\,sg\,\bigr)\;=\;r\mathbb{E}(f)\,+\,s\mathbb{E}(g)$
     \item If $f \;\geq\; g$, then             $\mathbb{E}(f)\,\geq\,\mathbb{E}(g)$ 
\end{enumerate}

Let $\dm$ be a strict linear preorder respected by $>$ on the set of all gambles. The expectation $\mathbb{E}$ is said to be $\dm$-\term{increasing} if $f \;\dm\; g$ implies
            $\mathbb{E}(f)\,>\,\mathbb{E}(g)$ for all  $f,g\in \mathscr{X}$.\qeddef{}
\end{definition}

  Analogues of de Finetti's fundamental results can be established, even if  all three requirements imposed by a comparative prevision system are dropped.
\begin{theorem} \label{bbEtheorem} Let $\mathscr{X}$ and $\mathscr{Y}$ be collections of gambles with $\mathscr{X}\subseteq\mathscr{Y}$, let $(\mathscr{L},\succsim)$ is be a complete preorder of a linear space $\mathscr{L}$ containing all constant gambles, and let $\dm$ be a strict linear preorder respected by $>$ on the set of all gambles.
  \begin{enumerate}[\thmlist \textsc{\textup{cs}}1, leftmargin=0.65in,labelsep=0.25in,itemindent=0.0in,itemsep=.15in,topsep=.1in]
        \item Any  $\dm$-coherent comparative expectation system  on $\mathscr{X}$ can be extended to one on $\mathscr{Y}\supseteq \mathscr{X}$.
        \item The comparative expectation system $\bigl(\mathscr{L},\succsim\bigr)$ is $\dm$-coherent if and only if it is a totally preordered linear space for which the asymmetric part $\succ$ of $\succsim$ respects $\dm$.
        \item The comparative expectation system
     $\bigl(\mathscr{L},\succsim\bigr)$ is $\dm$-coherent if and only if there is a totally ordered field extension $\Fset$ of the system real numbers and a $\dm$-increasing $\Fset$-valued expectation $\mathbb{E}$ on $\mathscr{L}$ that represents $\bigl(\mathscr{L},\succsim\bigr)$ in the sense that
     for all gambles $f,g\in\mathscr{L}$:\\
        \begin{center} $f\;\succsim\; g$ \quad if and only if \quad $\mathbb E(f)\;\geq\; \mathbb E(g)$.
        \end{center}
        \qedthm  
\end{enumerate}
\end{theorem}

Part \textsc{\thmlist{cs1}} is a corollary of an analogue of de Finetti's Fundamental Theorem of Prevision.  Subject to the additional assumption that all gambles from $\mathscr{L}$ are bounded,  a corollary of part \textsc{\thmlist{cs2}} is that  $\bigl(\mathscr{L},\succsim\bigr)$ is a comparative prevision structure if and only if it is an Archimedean totally preordered linear space for which the asymmetric part $\succ$ of $\succsim$ respects $\gg$.  As with comparative previsions, part \textsc{\thmlist{cs3}}  provides a numerical representation of comparative probability in terms of an $\Fset$-valued probability function when the $\Fset$-valued expectation is restricted to indicator functions for events. 


We conclude this section by clarifying for the reader how the present paper goes beyond \cite{Pedersen:2014}.  Blithely asserted in passing in that paper is, essentially, that part \textsc{\thmlist{cs3}} can be established by a routine ultraproduct construction and that the $\Fset$ may be taken to be a \emph{Hahn field of formal power series}.    The present paper not only gives rigorous definition to these claims but also substantiates them in the form of a newly obtained generalization of the Hahn Embedding Theorem.
Even excluding the applications and examples developed in the sequel, a significant contribution of the present paper is to make good on the assertions using methods presented in \cite{Pedersen:2014}.\footnote{\label{fn:specialissue}As such, the present paper is the product of substantial thinking beyond  \cite{Pedersen:2014}, which was prepared for a special issue honoring Horacio Arlo-Costa.} 

\section{Examples}
\label{sec:examples}

The next example provides simple illustrations of coherence and incoherence.

\begin{example}[Coherence and Incoherence] Consider binary relations $ \succsim_{1},\succsim_{2}$, and $\succsim_{3}$ over a common set of gambles $\mathscr{X}$ specified as follows.
\begin{myitemize}[leftmargin=3em,labelsep=1em]
\item[$\succsim_{1}$] $\varnothing\succsim_{1}\Omega$: Then  $(\mathscr{X},\succsim_{1})$ is $\dm$-incoherent.  To see this, observe that for $n=1$, gambles $g_{1}=\mathbf{1}$ and  $f_{1}=\mathbf{0}$ with  $c_{1}=1$ satisfy  $c_{1}g_{1} > c_{1}f_{1}$, whereby conditions \ccon[\dm]{1}, \ccon[\dm]{2}, and \ccon[\dm]{4} are satisfied.
\item[$\succsim_{2}$]  Either $E\sim_{2}\varnothing$ or $\Omega\sim_{2}E$ for some event $E\notin\{\Omega,\varnothing\}$: Then $(\mathscr{X},\succ_{2})$ is $>$-incoherent but both $\gg$-coherent and $\gtrdot$-coherent.
\item[$\succsim_{3}$] 
Both $\ind{H}-\ind{T}\succsim_{3} \mathbf 0$ and
$\ind{T}-\ind{H}\succsim_{3} \mathbf 0$ for pairwise disjoint events $H,T$ distinct from $\varnothing$.  Then $(\mathscr{X},\succ_{3})$ is $\dm$-coherent only if in addition both $\ind{H}-\ind{T}\sim_{3} \mathbf 0$ and $\ind{T}-\ind{H}\sim_{3} \mathbf 0$. 
To be sure, if, say, $\ind{H}-\ind{T}\succ_{3} \mathbf 0$, then
for $n=2$,  gambles $f_1=\ind{H}-\ind{T}=-f_{2}$ and $g_1=g_2=\mathbf 0$ with
$c_1=c_2=1$, 
conditions \ccon[\dm]{1}, \ccon[\dm]{2}, and \ccon[\dm]{3} are satisfied.\qedexa
\end{myitemize}
\end{example}

 The criterion of coherence in  Definition \ref{coherentComparativeExpectationDefn} can be generalized to parameterized criteria of $\Kset$-\emph{linear} $\dm$-\emph{coherence} for each totally ordered subfield $\Kset$ of $\Rset$ by stipulating that the coefficients $c_{1},\ldots,c_{n}$ additionally belong to $\Kset$ \citep{Pedersen:2014}.  The next example illustrates that the resulting criteria are not logically equivalent --- to wit, $\Qset$-linear $\dm$-coherence does not imply $\Rset$-linear
    $\dm$-coherence. All this notwithstanding,  the results of this paper generalize
    to these parameterized coherence criteria. 

\begin{example}[$\Kset$-linear Coherence]
    \label{hilberthotelexample}
    Let $\Omega$ be the set of positive rational numbers, and let $\alpha$ be a positive irrational number. Consider gambles $g_1$ and $g_2$  such that
     each positive rational number $q$:
\begin{center}
\begin{minipage}{0.45\textwidth}
\begin{align*}
          g_1(q) \quad&=\quad \begin{dcases}
            -\dfrac{q\,+\,\alpha }{2\alpha q} &\text{ if $q>\alpha$};\\[0.5em]
            \;\;\dfrac{q\,+\,\alpha }{2\alpha q} &\text{ if $q<\alpha$.}
        \end{dcases}
\end{align*}
\end{minipage}
\begin{minipage}{0.45\textwidth}
\begin{align*}
        g_2(q) \quad&=\quad \begin{dcases}
            \;\;1 &\mbox{ if $q>\alpha$};\\[0.5em]
            -1 &\mbox{ if $q<\alpha$.}
        \end{dcases}
\end{align*}
\end{minipage}
\end{center}
 \noindent Assume that $\boldsymbol{0}$ in addition to $g_{1},g_{2}$ belong to $\mathscr{X}$ and that $\succsim$ is fully determined by two judgments, viz., $\mathbf 0\succ g_1$
        and $\mathbf 0\succ g_2$.
    Then:
    \begin{myitemize}[leftmargin=0.65in,labelsep=0.25in,itemindent=0.0in,itemsep=.15in,topsep=.1in]

        \item[(i)] $(\mathscr{X},\succsim)$ is $\Qset$-linear $>$-coherent.
         \item[(ii)] $(\mathscr{X},\succsim)$ is $\Rset$-linear $\gtrdot$-incoherent; but
         \item[(iii)] $(\mathscr{X},\succsim)$ is $\Rset$-linear $\gg$-coherent.
    \end{myitemize}
 It follows from (i) that $(\mathscr{X},\succsim)$ is also both $\Qset$-linear $\gtrdot$-coherent and $\Qset$-linear $\gg$-coherent, while it follows from (ii) that $(\mathscr{X},\succsim)$  is also $>$-incoherent. We establish (i) and (ii).
\begin{myitemize}[leftmargin=2.5em,labelsep=1em]
\item[(i)] For \emph{reductio ad absurdum}, assume $\succsim$ is $\Qset$-linear incoherent.
        By cross-multiplying to eliminate denominators and simplifying, we may assume that  for $n=2$, gambles $f_1=f_2=\mathbf 0$ and $g_1$ $g_2$ defined as above with both positive $c_1,c_2\in \Qset$ satisfy conditions \ccon{1}, \ccon{2}, and \ccon{3$\mid$4}.  But $
            (\frac{c_{1}}{c_{2}})g_1(\nicefrac{c_{1}}{c_{2}})+g_2(\nicefrac{c_{1}}{c_{2}})< 0$ and so $c_{1}g_1(\nicefrac{c_{1}}{c_{2}})+c_{2}g_2(\nicefrac{c_{1}}{c_{2}})< c_{1}f_1(\nicefrac{c_{1}}{c_{2}})+c_{2}f_2(\nicefrac{c_{1}}{c_{2}})$, which conflicts with requirement \ccon{2}.
\item[(ii)]
For $n=2$, gambles $f_1=f_2=\mathbf 0$ and $g_{1}$ and $g_{2}$ defined as above with $c_1=\alpha$ and $c_2=1$ satisfy conditions \ccon[\gtrdot]{1}, \ccon[\gtrdot]{2}, and \ccon[\gtrdot]{3$\mid$4}.  To see that condition \ccon[\gtrdot]{4} and hence \ccon[\gtrdot]{2} obtain, observe that $\alpha g_1(q)+g_2(q)>\alpha\frac{\sgn{(\alpha -q)}}{\alpha}+\sgn{(q-\alpha )}=0$
for all positive $q\in\Qset$ (here $\sgn{(x)}$ denotes the sign of $x$).

        \qedexa
    \end{myitemize}
\end{example}

\begin{example}[Everywhere-Defined Coherent Extension of Real-Valued Expectation]
    \label{coherencyexample}
   Consider an $\Rset$-valued expectation $\mathbb{P}$ on a linear space $\mathscr {L}$ including all constant gambles. From Theorem \ref{bbEtheorem} it follows that $\mathbb{P}$ is $\dm$-increasing (e.g., a real-valued finitely expectation if $\gg$-increasing)  if and only if it can be can extended to an $\dm$-increasing $\Fset$-valued expectation $\mathbb{E}$ on the set of all gambles for some ordered field extension $\Fset$ of the field of real numbers. 
    \qedexa{}
\end{example}
Any total preorder of elementary events covering a finite state space can be represented by a real-valued probability function on the power set algebra.
The following example shows that such a representation may be obtained on any state space  by relaxing the requirement that probability function be real-valued.
\begin{example}
    \label{pigeonholeexample}
 Consider the family $\mathscr{X} = \Bigl\{{\{\omega\}}\,:\,\omega\in\Omega\Bigr\}$ of elementary events on a state space $\Omega$.  Suppose $(\mathscr{X},\succsim)$ is a complete preorder on $\mathscr{X}$.
     Then the comparative expectation system $(\mathscr{X},\succsim)$ is $>$-coherent.
    Thus, by Theorem \ref{bbEtheorem} it follows that 
    there is a $>$-increasing $\Fset$-valued probability function $\boldsymbol{p}$ on the power set algebra  $\mathscr{P}(\Omega)$ that represents $\succsim$ on $\mathscr{X}$ in the sense that 
    $\boldsymbol{p}(H_{1})\geq\boldsymbol{p}(H_{2})$ just in case
    $H_{1}\succsim H_{2}$ for all elementary events $H_{1}$ and $H_{2}$ in $\mathscr{X}$.

    To see that the comparative expectation system $(\mathscr{X},\succsim)$ is $>$-coherent, assume otherwise,
   whereby for some positive $n\in\Nset$ there are gambles $f_1,\ldots,f_n,
    g_1,\ldots,g_n\in\mathscr X$ of the form $f_i=\ind{\{a_i\}}$
and $g_i=\ind{\{b_i\}}$ and positive $c_1,\ldots,c_{n}\in \Rset$  satisfying conditions \ccon{1}, \ccon{2}, and \ccon{3$\mid$4}. We consider two cases provided for by condition \ccon{3$\mid$4}.
\begin{myitemize}[leftmargin=2em,labelsep=0.5em]
\item[\ccon{4}] Then there is $j=1,\ldots,n$ such that  $\displaystyle\sum_{i=1}^{n} c_i\ind{\{b_i\}}(b_{j})>\sum_{i=1}^{n} c_i\ind{\{a_i\}}(b_{j})$ , or else the left-hand side of inequality \ccon{4} would vanish.  For each $\omega\in \Omega$,
    let $A_\omega$ and $B_{\omega}$ denote the following sums:
 {   \begin{center}
\begin{minipage}{0.4\textwidth}
\begin{equation*}
          A_\omega \quad\coloneqq\quad \sum^{n}_{i=1}c_i\ind{\{a_{i}\}}(\omega)
\end{equation*}
\end{minipage}
\begin{minipage}{0.4\textwidth}
\begin{equation*}
          B_\omega \quad\coloneqq\quad \sum^{n}_{i=1}c_i\ind{\{b_{i}\}}(\omega) \qquad\tag{$\dagger$}
\end{equation*}
\end{minipage}
\end{center}
}
    Let $C\coloneqq\Bigl\{a_{i}\,:\,i\in\{1,\ldots,n\}\Bigr\}\cup\{b_{j}\}$. Observe:
    {
      \begin{center}
\begin{minipage}{0.3\textwidth}
\begin{equation*}
         \sum_{\omega\in C}A_{\omega} \quad=\quad \sum^{n}_{i=1}c_i
\end{equation*}
\end{minipage}
\begin{minipage}{0.3\textwidth}
\begin{equation*}
       \sum_{\omega\in C}B_{\omega} \quad\leq\quad \sum_{i=1}c_i
\end{equation*}
\end{minipage}
\end{center}
}
Since $B_{a_{i}}\geq A_{a_{i}}$  for each $i=1,\ldots n$ and $B_{b_{j}}> A_{b_{j}}$, it follows that:
   \begin{align*}
          \sum_{\omega \in C} B_{\omega}\quad&>\quad \sum_{\omega \in C} A_{\omega},
\end{align*}
which is impossible.

\item[\ccon{3}] We may assume that condition \ccon{4} fails to obtain. 
     By \ccon{3},
    there is $j= 1,\ldots, n$ such that
 $\{b_j\}\prec \{a_j\}$. We may assume that 
    $\{a_n\}\succsim\cdots\succsim \{a_1\}$ and
  that $\{a_j\}\succ \{a_{j-1}\}$ if $j\geq 2$ (if
   $\{a_j\}\sim \{a_{j-1}\}$, swap $a_j$ with $a_{j-1}$,
    $b_j$ with $b_{j-1}$, and $c_j$ with $c_{j-1}$).
     observe that no $\ell,m\in\{1,\ldots,n\}$ with $m\leq j\leq\ell$ are such that  $\{b_m\}\sim \{a_\ell\}$, for otherwise 
    $\{a_{\ell}\}\succsim \cdots\succsim \{a_m\}\succsim \{b_m\}\sim \{a_{\ell}\}$, which is inconsistent with
    $\{b_j\}\prec \{a_j\}$ (if $j=1$) or $\{a_j\}\succ \{a_{j-1}\}$
    (if $j\neq 1$).

    Now let the sums $A_{\omega}$ and $B_{\omega}$ be defined as before in ($\dagger$)   Then  each $\ell\geq j$:
 \begin{align*}
 A_{a_{\ell}}\quad&=\quad B_{a_{\ell}}\\
\quad&=\quad \sum^{n}_{i=j+1}c_i\ind{\{b_{i}\}}(a_{\ell})
\end{align*}
Now let $D\coloneqq\Bigl\{a_{i}\,:\,i\in\{j,\ldots,n\}\Bigr\}$.  Then \begin{align*}
        c_{j}+  \sum^{n}_{i=j+1} c_{i} \quad&\leq \quad \sum_{\omega \in D} A_{\omega}\\
        \quad&=\quad \sum_{\omega \in D} \sum^{n}_{i=j+1}c_i\ind{\{b_{i}\}}(\omega)\\
     \quad&\leq \quad     \sum^{n}_{i=j+1} c_{i}
\end{align*}
   But this is impossible since $c_j$ is positive.
    \qedexa
    \end{myitemize}
    
\end{example}

A real-valued probability analog of Example \ref{pigeonholeexample} does not obtain.
For example, $\Omega$ could be the cardinal $2^{2^{\aleph_0}}$
and $\succeq$ could be the usual ordinal number ordering.
Then Example \ref{pigeonholeexample} produces more distinct probabilities than there
are real numbers.

\section{Representability}
\label{sec:Representability}

Given a totally ordered set $\Gamma$, let $\Rset\big(\,{\Gamma}\,\bigr)$ denote the real linear space of real-valued functions with well-ordered support, along with the zero constant function:
 \begin{align*}
\Rset\big(\,{\Gamma}\,\bigr)\quad&\coloneqq\quad\vastl\{\,f\in  \Rset^{\Gamma}\,:\, \Bigl(\,f\neq 0\,\Bigr)\,\mbox{ is well-ordered in }\Gamma\, \vastr\}\,\cup\,\Bigl\{\,\mathbf{0}\,\Bigr\}.
\end{align*}
The set $\Bigl(\,f\neq 0\,\Bigr) = \Bigl\{\,\gamma\in\dom{f}:f(\gamma)\neq 0\,\Bigr\}$ is called the \term{support} of $f$. Endowed with pointwise-defined addition and scalar multiplication, the set $\Rset\big(\,{\Gamma}\,\bigr)$ becomes a real linear space with the constant function $\mathbf{0}\,=\,\delta\longmapsto 0:\Gamma\to\Rset$ serving as its additive identity. Additionally equip $\Rset\big(\,{\Gamma}\,\bigr)$ with a binary relation $\leqq$ such that $ f\;\leqq\; g $ if and only if:
\begin{align*}
\quad &\quad\quad f\;=\;g\;\quad\quad\textup{ \,or\, }\quad\quad f\Biggl(\,\min\Bigl(\,f\neq g\,\Bigr)\,\Biggr)\;<\;g\Biggl(\,\min\Bigl(\,f\neq g\,\Bigr)\,\Biggr).
 \end{align*}
 
 Then $\Rset\bigl(\,{\Gamma}\,\bigr)$ becomes a totally ordered linear space.  
 \begin{definition} The  \term{Hahn lexical function space} over totally ordered set $\Gamma$ is the totally ordered linear space $\Rset\bigl(\,\epsilon^{\Gamma}\,\bigr)$ on $\Rset\big(\,{\Gamma}\,\bigr)$ so equipped with addition and scalar multiplication and compatible total order.\qeddef{}
 \end{definition}
 
  Any given element $f$ from Hahn lexical function space $\Rset\bigl(\,\epsilon^{\Gamma}\,\bigr)$ may be written in the form:
  
\[
\sum_{\gamma\in\Gamma}f_{\gamma}\epsilon^{\gamma},
\]
where each coefficient $f_{\gamma}$ is the value $f(\gamma)$.  Addition and scalar multiplication are thereby the familiar operations performed on formal power series in unknown $\epsilon$. As the choice of notation suggests, index $\gamma$ precedes index $\gamma^{\prime}$ precisely when the inequality $\epsilon^{\gamma^{\prime}}\;\leqq\;\epsilon^{\gamma}$ obtains.

 Now suppose that $\Gamma$ is a totally ordered Abelian group. Additionally equip the Hahn lexical function space  $\Rset\big(\,\epsilon^{\Gamma}\,\bigr)$ with multiplication by way of convolution, whereby for every $f,g\in \Rset\big(\,{\Gamma}\,\bigr)$:
 \begin{align*}
     f\cdot g\quad&\coloneqq\quad \sum_{\gamma\in\Gamma}\sum_{\gamma=\alpha+\beta}f_{\alpha}g_{\beta}\epsilon^{\gamma}.
 \end{align*}
 The resulting system based on function space $\Rset\big(\,{\Gamma}\,\bigr)$ is a totally ordered field extension of the real number system under the natural mapping sending $0\longmapsto \mathbf{0}$ and each nonzero real number $r\longmapsto r\epsilon^{0}$.
 
   \begin{definition} The  \term{Hahn lexical field of power series} over totally ordered Abelian group $\Gamma$ is the totally ordered field  $\Rset\bigl(\bigl(\,\epsilon^{\Gamma}\,\bigr)\bigr)$ so based on $\Rset\big(\,{\Gamma}\,\bigr)$.\qeddef{}
 \end{definition}
 
 \begin{example}[Laurent Formal Power Series] The Hahn field $\Rset\bigl(\bigl(\,\epsilon^{\Zset}\,\bigr)\bigr)$ over the system of integers $\Zset$ is the classic field of \term{Laurent formal power series}.  Each number $f\in\Rset\bigl(\bigl(\,\epsilon^{\Zset}\,\bigr)\bigr)$ may be written in the form:
\[  \sum^{\infty}_{n=n_{0}}f_{n}\epsilon^{n},
 \]
 where $n_{0}\in\Zset$ and $f_{n}\in\Rset$ for all $n\geq n_{0}$.  \qedexa{}
 \end{example}
 
 \begin{example}[Levi--Civita Formal Power Series] The Hahn field $\Rset\bigl(\bigl(\,\epsilon^{\Qset}\,\bigr)\bigr)$ over the system of rational numbers is the classic field of \term{Levi--Civita formal power series}.  Each number $f\in \Rset\bigl(\bigl(\,\epsilon^{\Qset}\,\bigr)\bigr)$ may be written in the form:
\[
  \sum^{\infty}_{n=0}s_{r_{n}}\epsilon^{r_{n}},
\]
where $(\,r_{n}\,)_{n=0}^{\infty}$ is an unbounded strictly increasing sequence in the field of rational numbers $\Qset$.  Clearly  the field $\Rset\bigl(\bigl(\,\epsilon^{\Qset}\,\bigr)\bigr)$  is an ordered field extension (of an order-isomorphic copy) of the field of Laurent formal power series $\Rset\bigl(\bigl(\,\epsilon^{\Zset}\,\bigr)\bigr)$.\qedexa{}
\end{example}

\begin{definition}[Truncation] Let $\Gamma$ be a totally ordered set.  Given $\xi\in \Gamma$, define the \term{cut at index} $\xi$ to be the mapping $\mathbf{c}_{\xi}:\Rset\big(\,{\Gamma}\,\bigr)\to\Rset\big(\,{\Gamma}\,\bigr)$ such that for each $f\in \Rset\big(\,{\Gamma}\,\bigr)$ and $\delta\in\Gamma$:
 \begin{align*}
\Bigl[\,\mathbf{c}_{\gamma}f\,\Bigr](\delta)\quad&\coloneqq\quad\begin{cases}f(\delta)&\mbox{ if }\delta< \xi;\\
0& \textup{otherwise}.\\
\end{cases}
\end{align*}
\qeddef{}
\end{definition}
Thus the cut at index $\xi$ is a linear transformation that \emph{truncates} every function $f$ at $\xi$. For example, the cut at index $3$ of the Laurent formal power series \[  f=\sum^{\infty}_{n=0}n!\;\epsilon^{n},
 \]
 is $\mathbf{c}_{3}(f)\;=\;1\,+\,\epsilon\,+\,2\epsilon^{2}$.

 \begin{definition}[Archimedean Equivalence]  Let $\mathbb{V}$ be a totally ordered linear space. Define binary relations $\iless$ and $\bowtie$ on $\mathbb{V}$ by setting for all $u,v\in\mathbb{V}:$
\begin{align*}
u\;\iless \;v\qquad&\mbox{:if and only if}\qquad  n|{u}|\,<\,|v|\quad\mbox{for all}\;n\in\Nset. \\
u\;\bowtie \;v\qquad&\mbox{:if and only if}\qquad  u\not\iless v \mbox{ and }v\not\iless u. 
 \end{align*} 
 Vector $u$ is said to be \term{infinitesimal relative to} vector $v$, and vector $v$ \term{infinite relative to} vector $u$, when $u\;\bowtie \;v$. Vector $u$ is said to be \term{Archimedean equivalent} to vector $v$ if $u\;\bowtie \;v$.\qeddef{}
 \end{definition}

The binary relation $\bowtie$ is an equivalence relation on $\mathbb{V}\setminus\{\,0\,\}$. Let $\mathbb{V}_{\bowtie}$ denote the collection of $\bowtie$-equivalence classes of $\mathbb{V}\setminus\{\,0\,\}$.  The $\bowtie$-equivalence class to which an nonzero element  of $\mathbb{V}$ belongs is called its 
 \term{Archimedean class}.

 \begin{definition} Let $\mathbb{V}$ be a totally ordered linear space. Define a binary relation $\trianglelefteq$ on $\mathbb{V}_{\bowtie}$ by setting for all $U,V\in\mathbb{V}_{\bowtie}$:
 \begin{align*}
U\;\trianglelefteq \;V\qquad&\mbox{:if and only if}\qquad  U\;=\;V \qquad\mbox{or}\qquad v\,\iless\,u\quad\mbox{for some}\; u\in U\;\mbox{and}\;v\in V.
 \end{align*}
  \qeddef{}
 \end{definition}

 It is readily verified that the binary relation $\trianglelefteq$ is well-defined.
 \begin{definition}\label{def:selc} Let $\mathbb{V}$ be a totally ordered linear space. A subset $\mathfrak{A}$ of $\mathbb{V}\setminus\{0\}$ is said to be a \term{selection of Archimedean representatives} from $\mathbb{V}$ if it satisfies the following properties:
 \begin{myitemize}[leftmargin=0.65in,labelsep=0.25in,itemindent=0.0in,itemsep=.05in,topsep=.1in]
\item[{\textcolor[RGB]{119,30,17}{\textsc{ar}1}}]  For each nonzero $u\in V$, there is some $u'\in \mathfrak{A}$ such that  $u\bowtie u'$; and
\item[{\textcolor[RGB]{119,30,17}{\textsc{ar}2}}]  For all nonzero $u,u'\in \mathfrak{A}$, if $u\bowtie u'$, then $u = u'$.\qeddef{}
\end{myitemize}
 \end{definition}

To state the next result, we recall some terminology and introduce a last bit of notation. Given a subfield $\Kset$ of a field $\Kset'$, an element ${t}\in \Kset'$ is said to be \term{algebraic over} $\Kset$ if there is a polynomial $p(x)\;=\;k_{0}+k_{1}x+k_{2}x^2 +\cdots+ k_{n}x^{n}$ with $(k_{i})_{i=1}^{n}\in \Kset^{n}$ not all zero such that $p(t)\;=\;0$; the field $\Kset'$ is said to be a (proper) \term{algebraic} extension of $\Kset$ if every element of $\Kset'$ is algebraic over $\Kset$ (and $\Kset$ is a proper subset of $\Kset'$).  A field $\Kset$ is accordingly said to be \term{real-closed} if it can be made into an ordered field admitting no proper algebraic ordered-field extension, or as mentioned before, equivalently, (i) for every $t\in \Kset_{>0}$, $t=u^{2}$ for some $u\in \Kset$ and (ii) for each polynomial $p(x)\;=\;k_{0}+k_{1}x+k_{2}x^{2} +\cdots+ k_{n}x^{n}$ with coefficients in $\Kset$, $n$ odd, and $k_{n}$ nonzero, there is $t\in \Kset$ such that $p(t)=0$.  By the Artin--Schreier Theorem, every ordered field $\Kset$ has a unique real-closed ordered-field algebraic extension (up to field isomorphism), called its \term{real closure} (for more on real closures and real-closed fields, see, for example, the presentations of \citet[pp. 345-346]{ChangKeisler:1990}, \citet[pp. 451-457]{Lang:2002},  \citet[pp. 95-96]{Marker:2002}, or \citet[pp. 419-428]{Steinberg:2010}.

 Given a subring $\Sset$ of $\Rset$ and subset $U$ of linear space $\mathbb{V}$, let $\Sset[U]$ denote the $\Sset$-\term{span} of $U$:
\begin{align*}
\Sset\bigl[U\bigr] \quad&\coloneqq \quad \vastl\{\,v\in \mathbb{V}\,:\,v=\sum^{n}_{i=1}s_{i}w_{i}\;\mbox{for some positive integer}\,n\,\mbox{and}\, (s,w)\in \Sset^{n}\times U^{n}\,\vastr\}
     \end{align*}

 A totally order field $\Kset$  is said to be \term{Archimedean complete} if it admits no  proper extension $\Kset^{\prime}$ such that element $a'\in \Kset^{\prime}\setminus\Kset$ is Archimedean equivalent to some $a\in \Kset$ (cf. \citet[p. 55ff] {Fuchs:1963}).

\begin{theorem} \label{thm:HH}Let $\mathbb{V}$ be a totally ordered linear space, let $\mathfrak{A}$ be a selection of positive Archimedean representatives, and let $\mathfrak{u}\in \mathfrak{A}$.  There is a mapping $\psi:\mathbb{V}\to\Rset\Bigl(\Bigl(\,\epsilon^{\Zset[\mathfrak{u}-\mathfrak{A}]}\,\Bigr)\Bigr)$ satisfying the following properties$:$
\begin{enumerate}[\thmlist \textsc{hh}1,leftmargin=3em,labelsep=2em,itemindent=2em,itemsep=1em,topsep=.1in]
\item  $\psi\bigl(\,v\,\bigr)\;=\;\epsilon^{\mathfrak{u}-v}$ \quad for every $v\in \mathfrak{A} $ ;
\item $\mathbf{c}_{v}\Bigl(\img{\psi}\Bigr)\;\subseteq\;\img{\psi}$ \quad for every $v\in \Zset[\mathfrak{u}-\mathfrak{A}]$;
\item    $\Rset\Bigl(\Bigl(\,\epsilon^{\Qset[\mathfrak{u}-\mathfrak{A}]}\,\Bigr)\Bigr)$ is (isomorphic to) the real closure of $\Rset\Bigl(\Bigl(\,\epsilon^{\Zset[\mathfrak{u}-\mathfrak{A}]}\,\Bigr)\Bigr)$;

\item  For every $r\in \Rset$ and $v,w\in \mathbb{V}$:
\begin{align*}
\psi\bigl(\,r\,v\,+\,w\,\bigr)\quad&=\quad r\psi(v)\,+\,w;
\end{align*}
\item For every $v,w\in \mathbb{V}$:
\begin{align*}
v\,\leq\, w\qquad & \mbox{if and only if}\qquad \psi(v)\,\leqq\,\psi(w). 
\end{align*}
\noindent In addition, the field $\Rset\bigl(\bigl(\,\epsilon^{\Zset[\mathfrak{u}-\mathfrak{A}]}\,\bigr)\bigr)$ (the real-closed field $\Rset\bigl(\bigl(\,\epsilon^{\Zset[\mathfrak{u}-\mathfrak{A}]}\,\bigr)\bigr)$) is, up to order-isomorphism, the smallest (real-closed) Archimedean complete totally ordered field extension of $\Rset$ to include an order-isomorphic copy of the linear space $\mathbb{V}$.
\end{enumerate}
\qedthm{}
\end{theorem}

\proof Apply Hahn's Embedding Theorem \cite[Theorem 3.1]{HausnerWendel:1952} for real linear spaces to obtain a one-to-one, order-preserving linear function $F:\mathbb{V}\to\Rset\bigl(\,\epsilon^{\mathbb{V}_{\bowtie}}\,\bigr)$ such that $F(v)=\ind{\{v\}}$ and $\mathbf{c}_{v}\bigl(\img{F}\bigr)\;\subseteq\;\img{F}$ for each $v\in \mathfrak{A}$.  Observe that for every $v,w\in \mathfrak{A}$, the relation $\mathfrak{u}-v\leq \mathfrak{u}-w$ obtains just in case the relation $v/{\bowtie}\,\trianglelefteq\, w/{\bowtie}$ obtains, so $\Rset\bigl(\,\epsilon^{\mathbb{V}_{\bowtie}}\,\bigr)$ is linear order-isomorphic to $\Rset\bigl(\,\epsilon^{\mathfrak{u}-\mathfrak{A}}\,\bigr)$ under the mapping $G:\Rset\bigl(\,\epsilon^{\mathbb{V}_{\bowtie}}\,\bigr)\to\Rset\bigl(\,\epsilon^{\mathfrak{u}-\mathfrak{A}}\,\bigr)$ for which $G(f)(\mathfrak{u}-v)=f(v/{\bowtie})$ for each $v\in\mathfrak{A}$ and $f\in \Rset\bigl(\,\epsilon^{\mathbb{V}_{\bowtie}}\,\bigr)$. 

Now consider the Hahn lexical field  $\Rset\Bigl(\Bigl(\,\epsilon^{\Zset[\mathfrak{u}-\mathfrak{A}]}\,\Bigr)\Bigr)$ over the totally ordered $\Zset$-span $\Zset[\mathfrak{u}-\mathfrak{A}]$ of $\mathfrak{u}-\mathfrak{A}$ and the linear order-embedding  $H:\Rset\bigl(\,\epsilon^{\mathfrak{u}-\mathfrak{A}}\,\bigr)\to\Rset\bigl(\bigl(\,\epsilon^{\Zset[\mathfrak{u}-\mathfrak{A}]}\,\bigr)\bigr)$ such that for every $f\in \Rset\bigl(\,\epsilon^{\mathfrak{u}-\mathfrak{A}}\,\bigr)$ and $v\in \Zset[\mathfrak{u}-\mathfrak{A}]$:
 \begin{align*} 
   H(f)(v)\quad&\coloneqq\quad \begin{cases}f(v)&\textup{if }v\in \mathfrak{u}-\mathfrak{A};\\
   0&\textup{otherwise.}
   \end{cases}
   \end{align*}
Let $\psi\coloneqq H\circ G \circ F$. It is readily verified that $\psi$ is a one-to-one, order-preserving linear function of $\mathbb{V}$ into $\Rset\bigl(\bigl(\,\epsilon^{\Zset[\mathfrak{u}-\mathfrak{A}]}\,\bigr)\bigr)$ such that $\psi\bigl(\,v\,\bigr)\;=\;\epsilon^{\mathfrak{u}-v}$ for every $v\in \mathfrak{A} $ and $\mathbf{c}_{v}\bigl(\img{\psi}\bigr)\;\subseteq\;\img{\psi}$ for every $v\in \Zset[\mathfrak{u}-\mathfrak{A}]$. 

Because the totally ordered Abelian group $\Qset[\mathfrak{u}-\mathfrak{A}]$ is divisible, the Hahn lexical field $\Rset\bigl(\bigl(\,\epsilon^{\Qset[\mathfrak{u}-\mathfrak{A}]}\,\bigr)\bigr)$ is real closed.  Furthermore, since $\Qset[\mathfrak{u}-\mathfrak{A}]$ is (isomorphic to) the divisible closure (hull) of $\Zset[\mathfrak{u}-\mathfrak{A}]$, it follows that $\Rset\Bigl(\Bigl(\,\epsilon^{\Qset[\mathfrak{u}-\mathfrak{A}]}\,\Bigr)\Bigr)$ is (isomorphic to) the real closure of $\Rset\Bigl(\Bigl(\,\epsilon^{\Zset[\mathfrak{u}-\mathfrak{A}]}\,\Bigr)\Bigr)$ \citep[pp. 442-443, Theorem 5.2.12]{Steinberg:2010} (cf. presentation of \citet[p. 218, Observation 2]{Alling:1987}, \citet[pp. 354-355, Lemma 5.4.13]{ChangKeisler:1990}, \citet[p. 46, Theorem 2.15]{DalesWoodin:1996}.)

\endproof

\proof{Proof of Theorem \ref{bbEtheorem}.} Suppose that $(\mathscr{L},\succsim)$ is $\dm$-coherent.  Consider the equivalence relation $\sim$ on $\mathscr{L}$ such that $f\sim g$ if and only if both $f\succsim g$ and $g\succsim f$.  Let $\mathbb{V}_{\mathscr{L}}\coloneqq\mathscr{L}/\!\!\sim$ be the collection of $\sim$-equivalence classes $f/\!\!\sim$ endowed with addition $+_{\mathbb{V}_{\mathscr{L}}}$, scalar multiplication $\cdot_{\mathbb{V}_{\mathscr{L}}}$, and a total ordering $\succsim_{\mathbb{V}_{\mathscr{L}}}$ inherited from addition $+$, scalar multiplication $\cdot$\,, and $\succsim$ for $(\mathscr{L},\succsim)$.  Then $\mathbb{V}_{\mathscr{L}}$ is an ordered linear space that respects $\dm$.

  Let $\psi:\mathbb{V}_{\mathscr{L}}\to\Rset\bigl(\bigl(\,\epsilon^{\Zset[\mathfrak{u}-\mathfrak{A}]}\,\bigr)\bigr)$ be a linear embedding from Theorem \ref{thm:HH}, where $\mathfrak{u}\in \mathfrak{A}$ is such that $\mathfrak{u}\;=\;{1/\!\sim}$. Define $\mathbb{E}:\mathscr{L}\to \Rset\bigl(\bigl(\,\epsilon^{\Zset[\mathfrak{u}-\mathfrak{A}]}\,\bigr)\bigr)$ by setting $\mathbb{E}(f)\coloneqq\phi({f/\!\sim})$ for every $f\in \mathscr{L}$.  Then by construction it follows that $\mathbb{E}$ is a finitely additive expectation such that for every $f,g\in\mathscr{L}$:
\begin{eqnarray*}
f\;\succsim\; g\quad &\Longleftrightarrow&\quad f/\!\sim\;\;\;\succsim_{\mathbb{V}_{\mathcal{L}}}\;\; g/\!\sim\\[1em]
&\Longleftrightarrow&\quad\mathbb{E}(f)\;\geq\; \mathbb{E}(g).
\end{eqnarray*}

\endproof

\section{Final Remarks}
\label{sec:remarks}

Theorem \ref{thm:HH} adapts a proof of Hahn's Embedding Theorem for real linear spaces due to \citet[Theorem 3.1]{HausnerWendel:1952}. That proof proceeds roughly along the lines of the standard proof of the Hahn-Banach Extension Theorem, another fundamental result that \citet{Hahn:1927} proved (independently of \citet{Banach:1929a,Banach:1929b}).   The standard proof of the Hahn-Banach Extension Theorem proceeds by first demonstrating that a bounded linear functional defined on a subspace of a normed linear space can be extended to a bounded linear functional defined on the span of the collection obtained by adjoining a single vector to the subspace, thereupon applying a version of Zorn's Lemma to show that a bounded linear functional defined on a subspace of a normed linear space can be extended to a bounded linear functional defined on the entire linear space. \citet{Clifford:1954} has shown the basic approach \cite{HausnerWendel:1952} in fact readily generalizes to establish the general form of Hahn's Embedding Theorem for ordered Abelian groups. \citet[pp. 56-61]{Fuchs:1963} and \citet[pp. 53-60]{Alling:1987} offer proofs using valuations based on the proof of \citet{HausnerWendel:1952} and \citet{Clifford:1954}.

\citet{Hausner:1954} and \citet{Thrall:1954} applied and developed the results, concepts, and techniques of \citep{HausnerWendel:1952} to pursue an investigation of lexicographic expected utilities in the von Neumann-Morgenstern framework, although \citet[p. 631]{vonNeumannMorgenstern:1947} had briefly entertained the idea of lexicographic utility.  Thereafter, authors such as \citet[p. 209, pp. 215 ff.]{Chipman:1960} and \citet{Fishburn:1971, Fishburn:1982} investigated vector-valued ordinal and cardinal utilities.  \citet{Fishburn:1974} has provided a comprehensive mathematical survey of contemporary research in this vein, which primarily focuses on finite-dimensional vector-valued utilities.  
 
 The basic idea of lexicographic representations  of probability has been around for a while, mentioned by authors such as \citet[p. 41]{Savage:1954}, \citet[pp. 1184]{Fine:1971}, and \citet[pp. 1458-1460]{Fishburn:1974}, and lexicographic probability, and more generally lexicographic subjective expected utility in the Savage--Anscombe--Aumann framework, has been studied extensively in recent years by authors such as \citet{Segal:1986}, \citet{BlumeBrandenburgerDekel:1991a}, \citet{Hammond:1994, Hammond:1999}, \citet{Halpern:2010}, and \citet{brickhill2018triangulating}. Research on non-Archimedean representations of belief and value initiated around the 1950s has gradually grown over the years; \citet{BlumeBrandenburgerDekel:1989} and \citet{FishburnLavalle:1998} provide an overview of work on lexicographic utility and expected utility which includes discussion of more recent developments. \citet{BlumeBrandenburgerDekel:1989} also present their theory of finite-dimensional lexicographic expected utility \citep{BlumeBrandenburgerDekel:1991a}, while \citet{FishburnLavalle:1998} also present their theory of subjective expected utility combining vector-valued utilities with matrix-valued subjective probabilities.

  \citet{Narens:1974} has studied comparative probabilities with numerical representations constructed from ultraproducts, although \citet[p. 413]{KraftPrattSeidenberg:1959} and \citet[p. 45]{Richter:1971} expressed interest in this area in passing. More recent work by \citet[pp. 818-820]{ChuaquiMalitz:1983}, \citet[pp. 307-308]{Coletti:1990}, and \citet[pp. 211-213]{Domotor:1994} employ ultraproduct constructions to  establish representation results for comparative probabilities, although non-Archimedean representations are not the focus of these papers. \citet{Halpern:2010} focuses on mathematical relationships between lexicographic, conditional, and nonstandard probabilities, providing a survey of recent research in this area. 

\footnotesize{\section*{Acknowledgements}
In view of the long and winding road this paper has taken, the authors are indebted to many who have partaken in its history. The authors gratefully acknowledge Horacio Arló-Costa, Catrin Campbell-Moore, Clint Davis-Stober, Kenny Easwaran, Melissa Fusco, Haim Gaifman, Gerd Gigerenzer, Michael Grossberg, Ralph Hertwig, Frederik Herzberg, Leon Horsten, Thomas Icard, M. Ali Kahn, David Kellen, Jason Konek, Konstantinos Katsikopoulos, Amit Kothiyal, Hannes Leitgeb, Isaac Levi, Laura Martignon, Conor Mayo-Wilson, Enrique Miranda, Larry Moss, Louis Narens, Mel Nathanson, Jonathan Nelson, Eric Pacuit, Rohit Parikh, Richard Pettigrew, Marcus Pivato, Alex Pruss, Jan Willem-Romeijn, Dana Scott, Gerhard Schurz, David Schrittesser, Stanislav Speranski, Jack Stecher, Max Stinchcombe, Pat Suppes, Matthias Troffaes, Matthias Unterhuber, Peter Wakker, Sylvia Wenmackers, and Marco Zaffalon for their valuable feedback and insightful discussions drawing on results and ideas presented in drafts of this manuscript, spanning its early form prior to the interim hiatus through its subsequent resurrection.}


\end{document}